\documentclass{article}
\usepackage[utf8]{inputenc}
\usepackage{amsfonts,amsmath,amssymb,mathtools}
\usepackage[left=2cm,right=2cm,top=2cm,bottom=2cm,bindingoffset=0cm]{geometry}
\usepackage{graphicx}
\usepackage{authblk}
\usepackage{xcolor}
\usepackage{hyperref}
\usepackage{authblk}
\usepackage{slashed}
\usepackage{braket}
\usepackage{comment}
\usepackage{pb-diagram}
\usepackage{wasysym}
\usepackage{amsthm}
\usepackage{cancel}
\usepackage{tikz}\usetikzlibrary{tikzmark}
\usetikzlibrary{fadings} 
\usetikzlibrary{positioning}
\usetikzlibrary{backgrounds} 
\usetikzlibrary[backgrounds]
\usetikzlibrary{decorations.pathreplacing}
\usetikzlibrary{arrows}
\definecolor{airforceblue}{rgb}{0.36, 0.54, 0.66}	
\definecolor{beige}{rgb}{0.96, 0.96, 0.86}
\definecolor{bittersweet}{rgb}{1.0, 0.44, 0.37}
\definecolor{melon}{rgb}{0.99, 0.74, 0.71}
\definecolor{mustard}{rgb}{1.0, 0.86, 0.35}
\definecolor{lava}{rgb}{0.81, 0.06, 0.13}
\definecolor{magnolia}{rgb}{0.97, 0.96, 1.0}
\definecolor{lavendermist}{rgb}{0.9, 0.9, 0.98}
\definecolor{lavendergray}{rgb}{0.77, 0.76, 0.82}
\definecolor{palepink}{rgb}{0.98, 0.85, 0.87}
\definecolor{palesilver}{rgb}{0.79, 0.75, 0.73}
\definecolor{cadetgrey}{rgb}{0.57, 0.64, 0.69}
\definecolor{anti-flashwhite}{rgb}{0.95, 0.95, 0.96}
\colorlet{Light0anti-flashwhite}{anti-flashwhite!70!white}
\colorlet{Lightanti-flashwhite}{anti-flashwhite!50!white}
\colorlet{Light2anti-flashwhite}{anti-flashwhite!30!white}
\definecolor{linkcolor}{rgb}{0,0,1}
\definecolor{urlcolor}{rgb}{0,0,1}

\usepackage{tikz-cd}
\usetikzlibrary{decorations.markings}
\usetikzlibrary{positioning}
\usepackage{array}
\usepackage{ytableau}
\usepackage{longtable}
\usepackage{empheq}
\usepackage{multicol}
\usepackage{mathrsfs}
\usepackage{diagbox}
\usepackage[vcentermath]{youngtab}
\usepackage[natbib=true, backend = bibtex, style=numeric-comp, sorting=none]{biblatex}
\hypersetup{pdfstartview=FitH,  linkcolor=linkcolor,urlcolor=urlcolor,citecolor=urlcolor, colorlinks=true}

\newcommand\bem{\begin{pmatrix}}
\newcommand\eem{\end{pmatrix}}
\newcommand\beq{\begin{equation}}
\newcommand\eeq{\end{equation}}
\newcommand\beqs{\begin{equation*}}
\newcommand\eeqs{\end{equation*}}

\date{}
\begin{document}

\title{\bf On an alternative stratification of knots}
\author[1,2,4]{{\bf E.~Lanina}\thanks{\href{mailto:lanina.en@phystech.edu}{lanina.en@phystech.edu}}}
\author[1,2,3,4]{{\bf A.~Popolitov}\thanks{\href{mailto:sleptsov@itep.ru}{popolit@gmail.com}}}
\author[1,2,4]{{\bf N.~Tselousov}\thanks{\href{mailto:tselousov.ns@phystech.edu}{tselousov.ns@phystech.edu}}}

\vspace{5cm}

\affil[1]{Moscow Institute of Physics and Technology, 141700, Dolgoprudny, Russia}
\affil[2]{NRC "Kurchatov Institute", 117218, Moscow, Russia}
\affil[3]{Institute for Information Transmission Problems, 127994, Moscow, Russia}
\affil[4]{Institute for Theoretical and Experimental Physics, 117218, Moscow, Russia}
\renewcommand\Affilfont{\itshape\small}
\maketitle

\vspace{-7cm}

\begin{center}
	\hfill ITEP-TH-18/22\\
	\hfill IITP-TH-17/22\\
	\hfill MIPT-TH-15/22
\end{center}

\vspace{5cm}

\begin{abstract}
We introduce an alternative stratification of knots: by the size of lattice on which a knot can be first met. Using this classification, we find ratio of unknots and knots with more than 10 minimal crossings inside different lattices and answer the question which knots can be realized inside $3\times 3$ and $5\times 5$ lattices. In accordance with previous research, the ratio of unknots decreases exponentially with the growth of the lattice size. Our computational results are approved with theoretical estimates for amounts of knots with fixed crossing number lying inside lattices of given size.
\end{abstract}

\setcounter{equation}{0}
\section{Introduction}\label{intro}
The study of knots representations is a central theme in knot theory. We are especially interested in the so called universal representations, which contain all knots. These ones include minimum braid representations, Morse link representations, arc representations, and many others. However, for example, arborescent knots are not universal, as an arbitrary knot cannot be represented in an arborescent form.

Different knots representations are useful for different purposes. For example, braid representation allows one to calculate quantum knot invariants using Reshetikhin-Turaev algorithm~\cite{guadagnini1990chern,reshetikhin1990ribbon,Morozov_2010,Kaul_1992,Rama_Devi_1993,Ramadevi_1994,Ramadevi_2001,arxiv.1107.3918,arxiv.1209.1346,Nawata_2014,Gu_2015,deguchi2014exchange,Mironov_2012,ITOYAMA_2012,ITOYAMA_2013,Anokhina_2013,anokhina2014cabling}. The HOMFLY polynomials for arborescent knots are more easier to calculate in terms of modular transformation matrices $S$ and $T$ and their conjugates~\cite{2RAMADEVI_1994,nawata2013colored,galakhov2015colored,galakhov2015knot,mironov2015colored,mironov2015colored}. 

In this paper, we investigate classical knot theoretic questions with the use of the so called lattice representation. A lattice representation is defined on a square lattice of size $2n+1$. Each node of the lattice is equipped with one of two crossings which we denote as $+$ and $-$. It was proved that each knot has at least one lattice diagram~\cite{Zohar}, i.e. this representation is universal. Moreover, each lattice knot in a fixed size lattice also can be realized inside any bigger lattice. So, we can introduce knots representation by lattice diagrams (see Section~\ref{defs} also). 
\begin{figure}[!ht]
    \centering
    \includegraphics[scale=0.29]{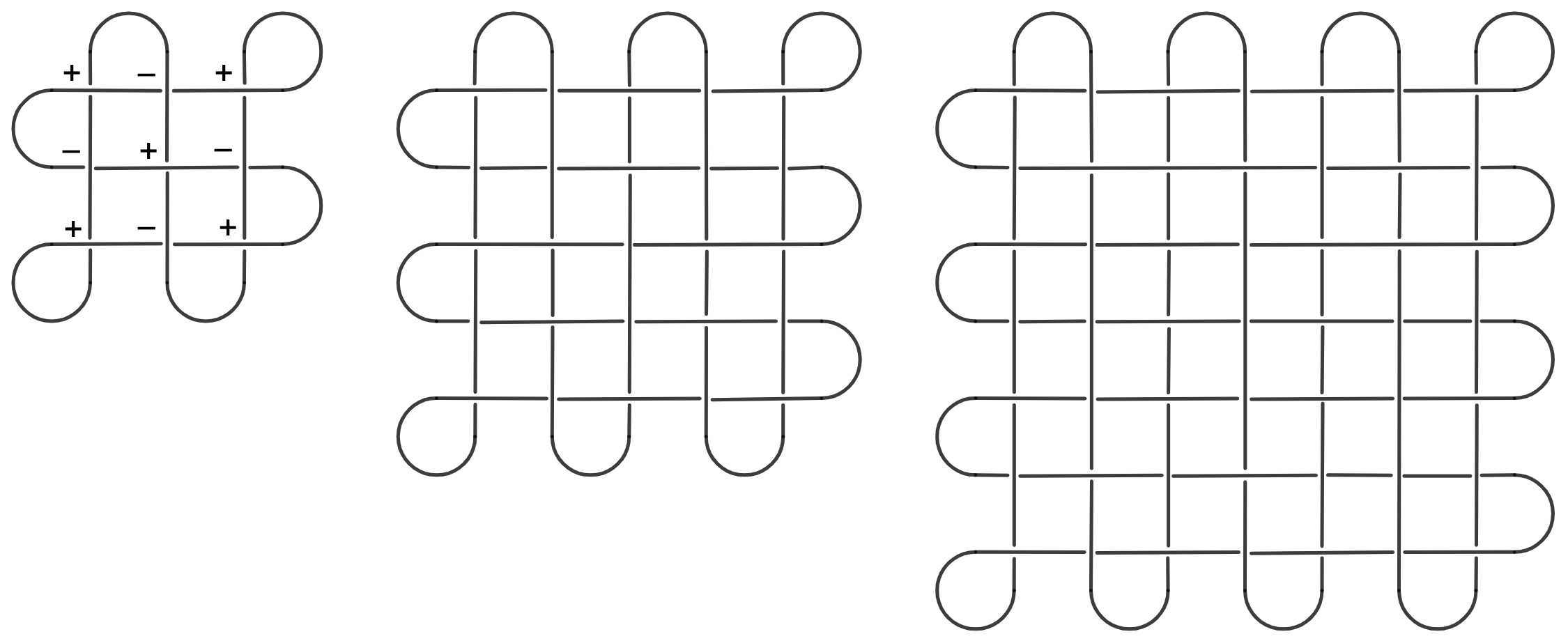}
\caption{Examples of lattice knots: $3 \times 3$, $5 \times 5$ and $7 \times 7$ lattice knots.}\label{LK}
\end{figure}\\

The set of all distinct knots $\mathbb{K}$ is infinite and to deal with this infinity the set is usually
stratified (split) into finite parts $\mathbb{K}_n,\ n = 1 \dots \infty$ according to the value
of some simple knot invariant. Commonly used is the splitting by the crossing number, i.e. the minimal number of crossings among all the knot representations for a given knot. Knots with small crossing numbers are collected in the celebrated Rolfsen table (see, for example,~\cite{katlas}). 

After a stratification of the knot set $\mathbb{K}$ is chosen, one can ask how
this splitting interplays with some other splitting of the knot set into different classes, for instance,
what is an asymptotic of a ratio of different knot classes (w.r.t the latter splitting)
in the $\mathbb{K}_n$ stratum as $n \rightarrow \infty$.

Indeed, following famous theorem of W. Thurston, all knots are divided into three types: torus, satellite and hyperbolic. For knots classification by crossing number $n$, it was proven that hyperbolic knots do not dominate for $n\rightarrow\infty$, and it was argued that satellite knots dominate~\cite{belousov2019hyperbolic}. However, this result is rather unexpected because based on numeric evidence for small $n$ for a long time it was believed that hyperbolic knots dominate for infinitely big crossing numbers due to a well-known conjecture (see~\cite{adams1995knot}, p.119).
It was also obtained that the number of prime knots grows exponentially in $n$~\cite{ernst1987growth}.
Given an alternative stratification of knots, one can ask whether these counterintuitive statements
are stable/sensitive w.r.t. change of stratification. In particular, which type of knots will dominate, if we split knots by the minimal size of their respective lattice diagram? What will be the distribution of knots coming from a fixed size lattice? In this short note we start programs to answer these questions (see Section~\ref{NumComp}) and obtain some theoretical bounds on numbers of different types of knots in a fixed size lattice (see Section~\ref{Bounds}).

Namely, main results of this paper are the following ones. First, we derive upper and lower bounds for the numbers of knots with fixed crossing numbers which can be realized inside a lattice of fixed size (Section~\ref{Bounds}). Second, we approve our bounds by numerical computations of knots that can be obtained from $3\times 3$ and $5 \times 5$ lattice diagram (Section~\ref{3x35x5}). In addition, we calculate the number of unknots and knots with crossing numbers greater than 10 inside bigger lattices (Section~\ref{sec:ratio-of-unknots}).

Lattice diagrams are also interesting because they connect knot theory and statistical models. The relation between knot theory and statistical mechanics was first noted in the work of Vaughan Jones~\cite{vf1985polynomial}. In his work, a connection between the Jones polynomials and the Potts model was established. Later, Jones developed a method to compute the HOMFLY polynomials using vertex models. This method was also used by Turaev for the Kauffman polynomials~\cite{turaev1990yang}. The connection between knot theory and statistical mechanics was formalised and further extended by Jones~\cite{jones1989knot} with the use of spin models.

The connection between exactly solvable statistical models and knot theory is interesting since it can lead to progress both in proving mathematical theorems in knot theory (for example, in the question of the dominance of satellite knots~\cite{belousov2019hyperbolic}) and in a search for new knot invariants and simpler methods for computing well-known invariants. We hope to move forward in this direction since the theory of integrable models of statistical physics is well developed~\cite{baxter1980exactly, lieb1972two}.


\setcounter{equation}{0}
\section{Basic theorems}\label{defs}
In this section, we introduce some basic facts on lattice diagrams to rely on in what follows. One can see recent paper~\cite{Zohar} for detailed explanations.

\textbf{Definition 1.} Lattice diagram $L_n$ is a closed immersed curve which has $2n+1$ horizontal and $2n+1$ vertical line segments and $(2n+1)^2$ crossing points.

In~\cite{Zohar}, lattice diagrams are called potholder curves or potholders. Knots carried by potholder curves were studied by Grosberg and Nechaev~\cite{Grosberg_1992,nechaev1999statistics}, who calculated the number of unknots carried by such curves via a connection to the Potts model of statistical mechanics.

A lattice knot is obtained from a lattice diagram by resolving each crossing to be one of two types which we denote as $+$ and $-\,$. The examples of lattice knots coming from $L_1$, $L_2$ and $L_3$ are shown in Fig.~\ref{LK}.

\textbf{Theorem 1.} Every knot is carried by lattice knots.

\textbf{Theorem 2.} All lattice knots coming from $L_n$ also come from $L_{n+1}$.

From these two theorems, we conclude that there is another classification of knots. Namely, we can sort knots by minimal size of lattice diagrams which they are produced by.


\setcounter{equation}{0}
\section{How can one distinguish knots?}
The main question we are interested in is which types of knots come from lattice diagram $L_n$. For this purpose, we discuss how one can distinguish lattice knots.

\subsection{Reidemeister moves}
The simplest, but laborious way is to use Reidemeister theorem.

\begin{figure}[h]
\begin{minipage}[h]{0.32\linewidth}
\center{\includegraphics[width=0.94\linewidth]{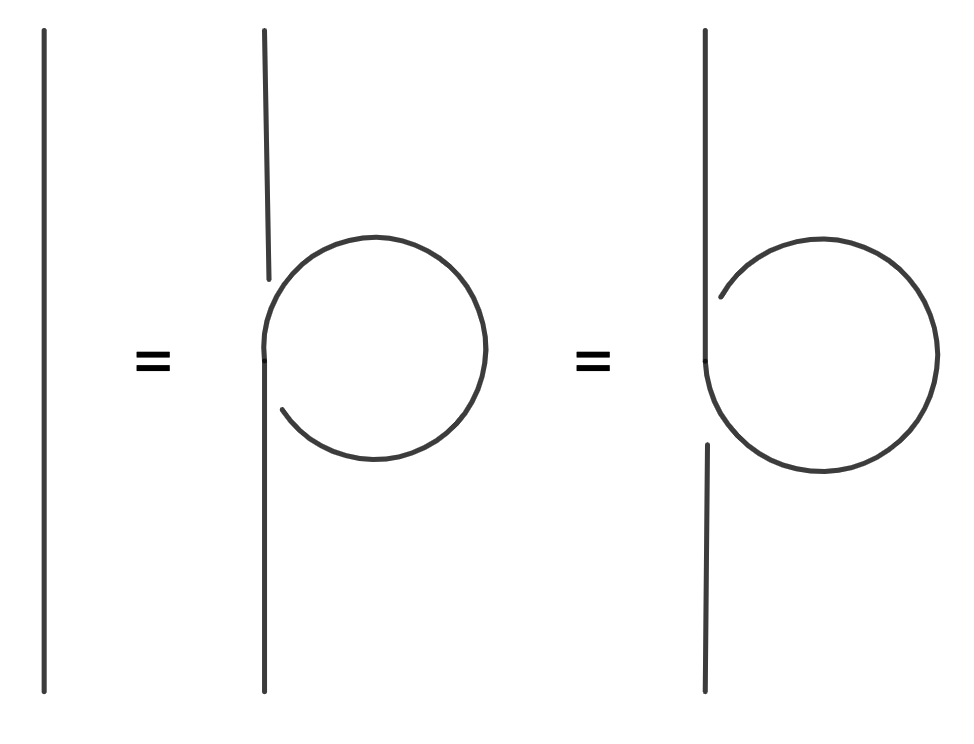}}
\end{minipage}
\hfill
\begin{minipage}[h]{0.32\linewidth}
\center{\includegraphics[width=0.94\linewidth]{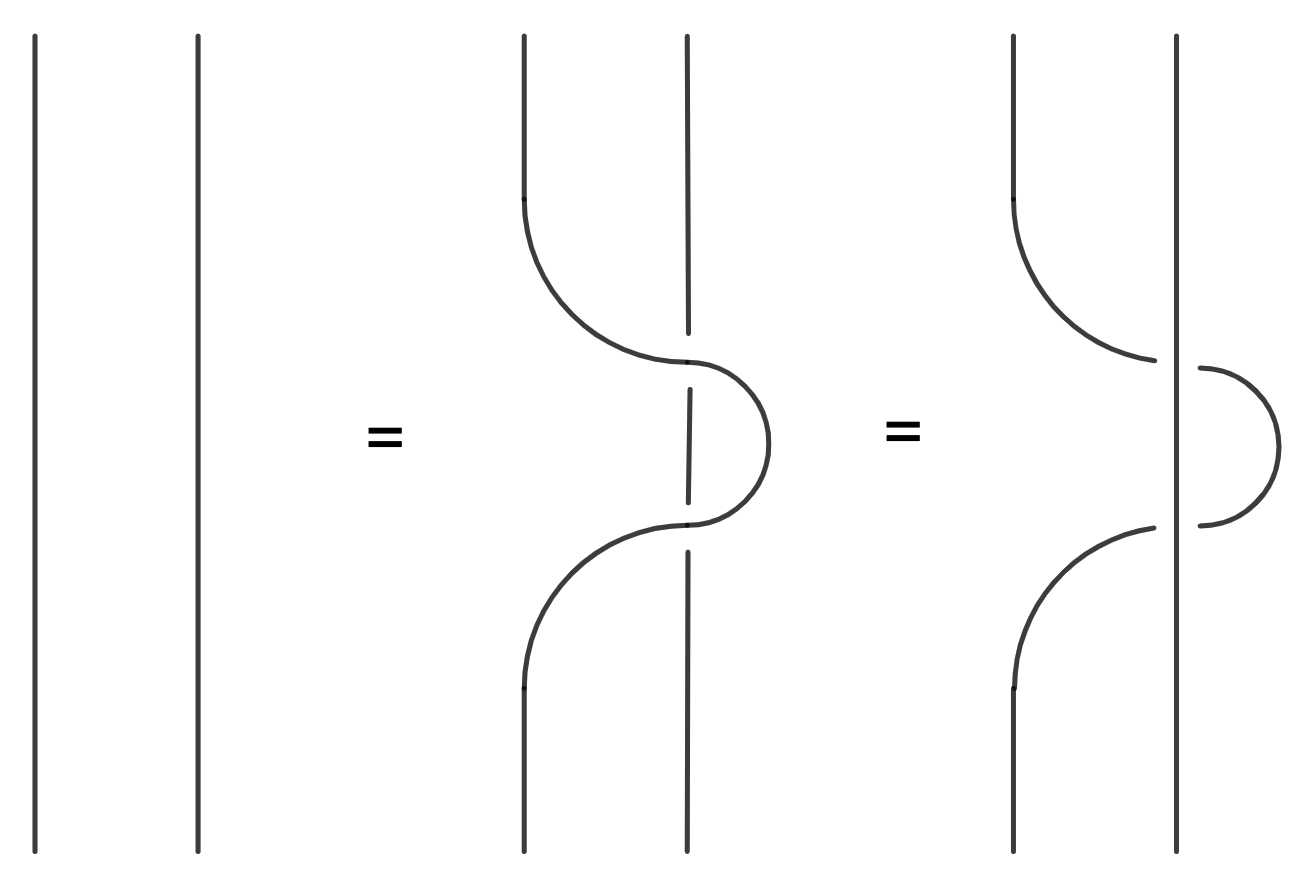}}
\end{minipage}
\hfill
\begin{minipage}[h]{0.32\linewidth}
\center{\includegraphics[width=0.94\linewidth]{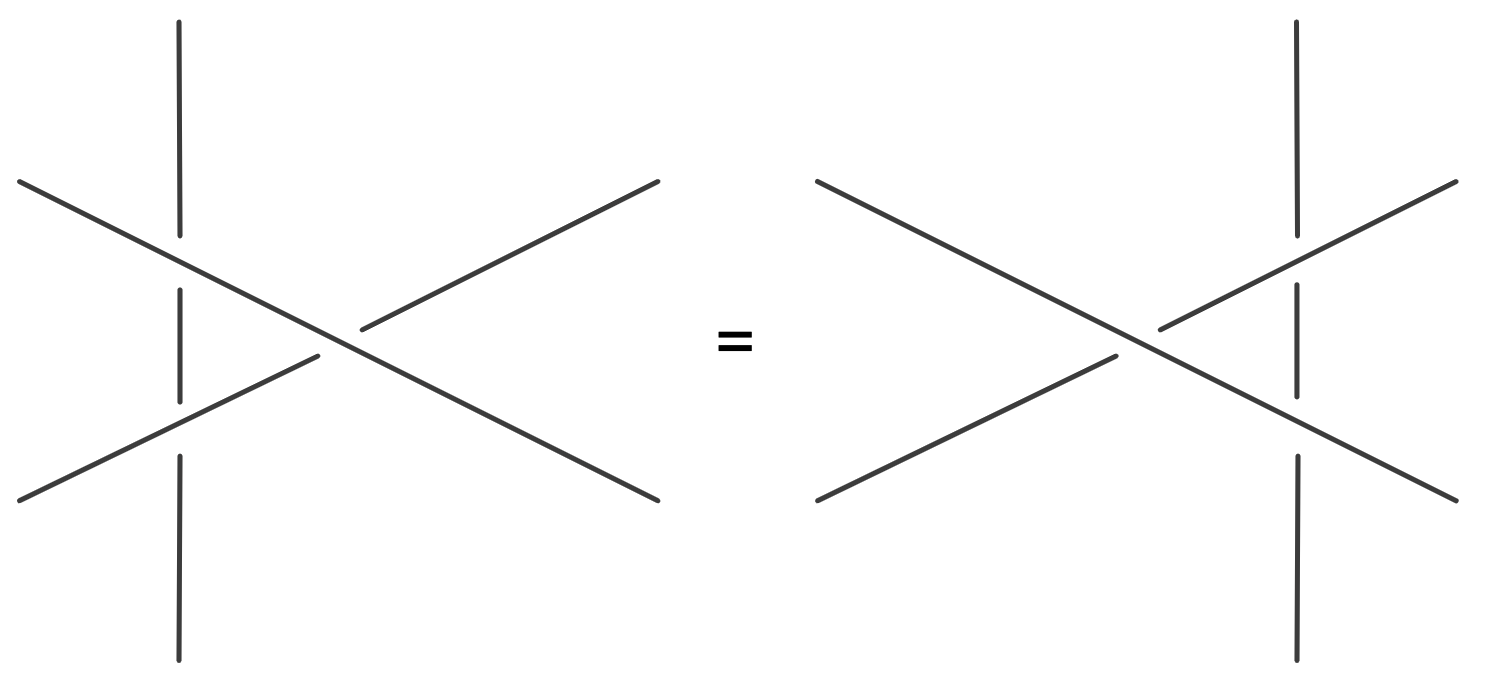}}
\end{minipage}
\caption{Reidemeister moves}\label{RM}
\end{figure}

\textbf{Theorem 3.} Two knots $\mathcal{K}_1$ and $\mathcal{K}_2$ are equivalent if there exist a sequence of Reidemeister moves (Fig.~\ref{RM}) which turn one projection into another.
\\
This theorem allows one to obtain from a lattice knot the same knots but inside lattices of bigger size. We use this approach to present theoretical bounds (Section~\ref{Bounds}).

\subsection{Jones polynomial}\label{Jones}
Knots can be practically distinguished with the help of peculiar special functions -- {\it knot invariants}. 

\textbf{Definition 2.} A knot invariant $I(\mathcal{K})$ is a function of a knot $\mathcal{K}$ which takes the same value on equivalent knots:
\begin{equation}
    \mathcal{K}_1=\mathcal{K}_2 \Rightarrow I(\mathcal{K}_1)=I(\mathcal{K}_2)\,.
\end{equation}
In our investigation, we utilize the following property of knot invariants:
\begin{equation}
    \mathcal{K}_1\neq\mathcal{K}_2 \Leftrightarrow I(\mathcal{K}_1)\neq I(\mathcal{K}_2)\,.
\end{equation}
There are many known knot invariants of different types in literature, we are interested in polynomial ones, because they can be effectively computed with the help of special computer programs. One of the most famous polynomial knot invariants is the Jones polynomial $J^{\mathcal{K}}$:
\begin{equation}
    J^{\mathcal{K}} : \mathcal{K} \rightarrow \mathbb{Z}[q, q^{-1}]
\end{equation}
that for a knot $\mathcal{K}$ returns a Laurent polynomial $J^{\mathcal{K}}(q)$. We provide examples of the Jones polynomials for several knots with small crossing number:
\begin{align}
    \begin{aligned}
    J^{0_1} &= 1\,, \\
    J^{3_1} &= q^2 + q^6 - q^8\,, \\
    J^{4_1} &= 1 -q^2 - q^{-2} + q^4 + q^{-4}\,.
    \end{aligned}
\end{align}
We consider the {\it reduced} or {\it normalized} Jones polynomial, therefore its value for unknot is $1$ rather than $q + q^{-1}$.
The Jones polynomial is a very good tool to distinguish knots with small number of crossings. It distinguishes almost knots in the Rolfsen table up to 9 crossings. The first example of knots which cannot be distinguished with the use of the Jones polynomial is
\begin{equation}
    J^{5_1}(q)=J^{10_{132}}(q)\,.
\end{equation}
Jones polynomial possesses a number of peculiar properties:
\begin{itemize}
    \item at the point $q = 1$ the Jones polynomial of any knot reduces to one:
    \begin{equation}
       J^{\mathcal{K}} (q = 1)  = 1\,.
    \end{equation}
    \item the Jones polynomial of the mirror image of a knot $\mathcal{K}^{mir}$ is expressed by the following formula through the Jones polynomial of the initial knot $\mathcal{K}$:
    \begin{equation}\label{JonesMir}
        J^{\mathcal{K}^{mir}}(q) = J^{\mathcal{K}}(q^{-1})\,.
    \end{equation}
    \item the Jones polynomial of a composite knot $\mathcal{K}_1 \# \mathcal{K}_2$ factorizes:
    \begin{equation}\label{JonesFactor}
        J^{\mathcal{K}_1 \# \mathcal{K}_2}(q) = J^{\mathcal{K}_1}(q) \cdot J^{\mathcal{K}_2}(q)\,. 
    \end{equation}
    \item the Jones polynomial of a disjoint union of two knots $\mathcal{K}_1 \bigsqcup \mathcal{K}_2$:
    \begin{equation}
        J^{\mathcal{K}_1 \bigsqcup \mathcal{K}_2}(q) = (q + q^{-1})\cdot J^{\mathcal{K}_1}(q) \cdot J^{\mathcal{K}_2}(q)
    \end{equation}
\end{itemize}
We get use of properties~\eqref{JonesFactor} and~\eqref{JonesMir} extensively in our computations.

\medskip

The Jones polynomial has several equivalent definitions and 
approaches to its computation. In our computations, we use {\it state sum formula} or {\it Khovanov algorithm}~\cite{Kh1,Bar_Natan_2002}, that is connected with rapidly growing area of knot homology theory. 

Here, we briefly present this algorithm to state its usefulness for distinguishing lattice knots. First, given a knot diagram, one writes down all knot smoothings (each intersection can be substituted with either 0-, or 1-smoothing); they become vertices of the so called Khovanov cube. Each knot smoothing correspond to a summand $(-1)^{r} q^{r}\left(q+q^{-1}\right)^{k}$, where $r$ is a number of 1-smoothings, while $k$ is a number of components of a knot smoothing. Then, one sums contributions of each knot smoothings and multiply the resulting polynomial by $(-1)^{n_{-}}q^{n_{+}-2n_{-}}(q+q^{-1})^{-1}$, where $n_-$ is a number of '$-$' crossings and $n_+$ is a number of '$+$' crossings. 

\begin{figure}[!ht]
    \centering
    \includegraphics[scale=0.52]{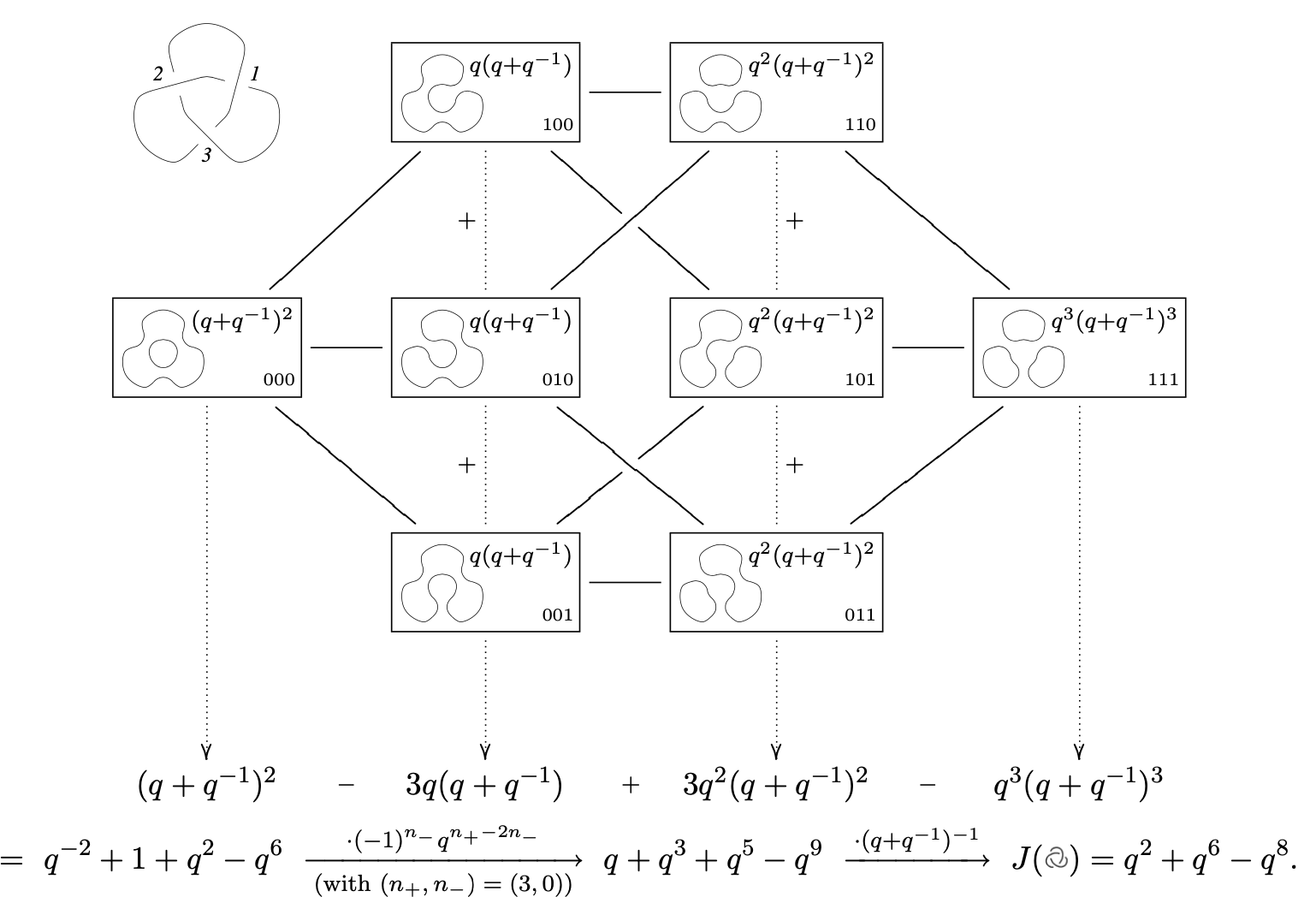}
\caption{Khovanov algorithm for calculation of the Jones polynmial for the trefoil from~\cite{Bar_Natan_2002}.}\label{trefoil}
\end{figure}

This algorithm becomes clear for the example of the trefoil (Fig.~\ref{trefoil}).
\\
\textbf{Remark 1.} Khovanov cubes are the same for all lattice knots coming from a lattice diagram of fixed size.

This fact allows us to effectively calculate the Jones polynomials for lattice knots. Namely, one needs to construct Khovanov cube for a fixed size lattice only once, and then, in order to obtain the Jones polynomials for all lattice knots inside a lattice with fixed size, one starts Khovanov algorithms with each of cube's vertex. That is exactly how our computer program from Section~\ref{3x35x5} works.

\setcounter{equation}{0}
\section{Upper and lower bounds}\label{Bounds}
For any $n$ and $k$, $n>k$, the types of knots realized on a lattice of size $(2k+1)\times(2k+1)$ can also be obtained on a lattice of size $(2n+1)\times(2n+1)$. Thus, we introduce another enumeration of knots analogous to the Rolfsen table. 

Namely, each knot is enumerated by two numbers: $m$ and $l_m$. The number $m$ corresponds to a minimal lattice size $(2m+1)\times(2m+1)$, on which the knot can be met. $l_m$ just enumerates distinct knots which can be represented as lattice knots with fixed minimal $m$. For example, for the trivial knot, we have $m=0$, $l_0=1$, and for the trefoil, we have $m=1$, $l_1=1$.

We want to know which knot dominates (i.e. has more representations) on a fixed size lattice. In a recent paper~\cite{belousov2019hyperbolic}, it was proved that the hyperbolic knots with a crossing number $\leq n$ do not dominate as $n\rightarrow \infty$ for the standard stratification of knots, and it was argued that satellite knots dominate instead. Will the picture be different if one classifies knots by lattices sizes?

In this section, we evaluate how many different knots can be met inside a fixed size lattice. Namely, we establish the connection between the introduced stratification by lattice diagrams and by the Rolfsen table. In other words, we answer the question, how many knots $n_{m}$ given from the Rolfsen table come from a lattice diagram $L_k$.

In order to evaluate numbers of different lattice knots, we split the problem into two parts:
\begin{enumerate}
    \item First, in Subsection~\ref{LatEmb}, we evaluate how many times one can find the $(2k+1)\times(2k+1)$ lattice inside the $(2n+1)\times(2n+1)$ lattice for $n>k$.
    \item Second, in Subsection~\ref{KnotsLat}, we define which knots classified by the Rolfsen table can be found in the $(2k+1)\times(2k+1)$ lattice.
\end{enumerate}
This lets us factor out the trivial combinatorial contribution and count only essentially distinct representations.

\subsection{Lattices embedding}\label{LatEmb}
In order to figure out how many ways we can embed a fixed lattice $(2k+1)\times(2k+1)$ into a lattice of size $(2n+1)\times(2n+1)$, we 'untie'\, the bigger lattice up to a lattice of a smaller size.

\begin{figure}[!ht]
    \centering
    \includegraphics[scale=0.35]{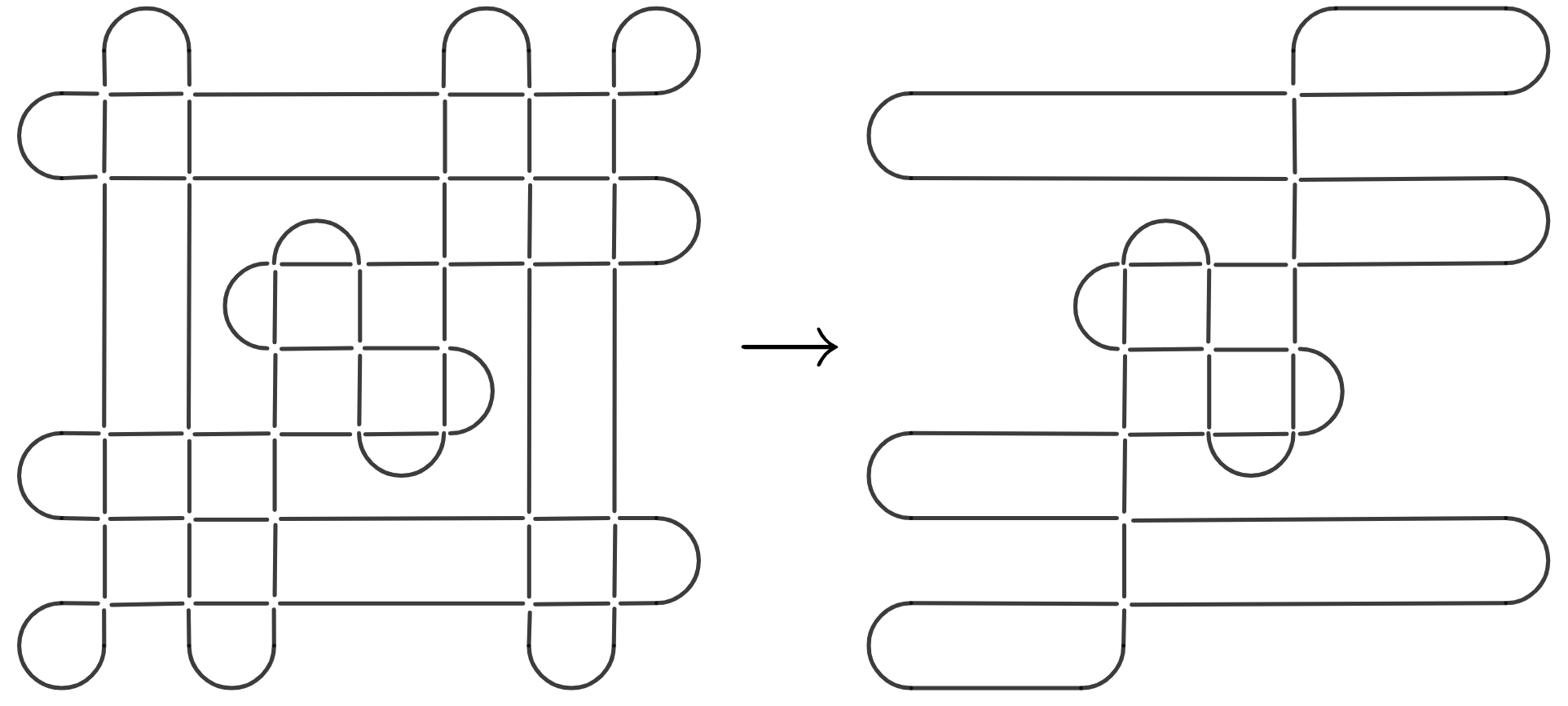}
    \caption{Smaller lattice knot in the middle of the bigger lattice.}\label{NotOnBoundary}
\end{figure}

Note that we can pull off a loop which is intersected by $m$ threads within $2^{m}$ options so that on each thread, we have the same type of intersections. Here are some examples:
\begin{figure}[!ht]
    \centering
    \includegraphics[scale=0.3]{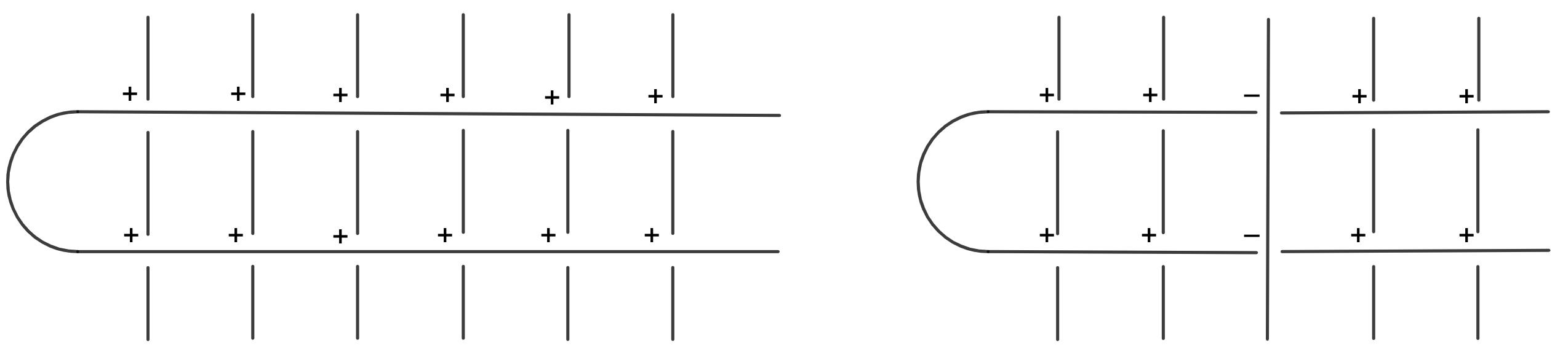}
\caption{Examples of loops pull off}
\end{figure}\\
Thus, to obtain the $(2k+1)\times(2k+1)$ lattice from the $(2n+1)\times(2n+1)$ lattice, we follow three steps.
\begin{enumerate}
    \item Pull off $k$ horizontal and $k$ vertical loops which belong to the $(2k+1)\times(2k+1)$-lattice. We can do this in $2^{2(n-k)\cdot 2k}$ ways.
    \item Then, we pull off $n-k$ vertical loops in $2^{2(n-k+1)(n-k)}$ variants. We provide the following examples:
    \begin{figure}[h]
\begin{minipage}[h]{0.5\linewidth}
\center{\includegraphics[width=0.9\linewidth]{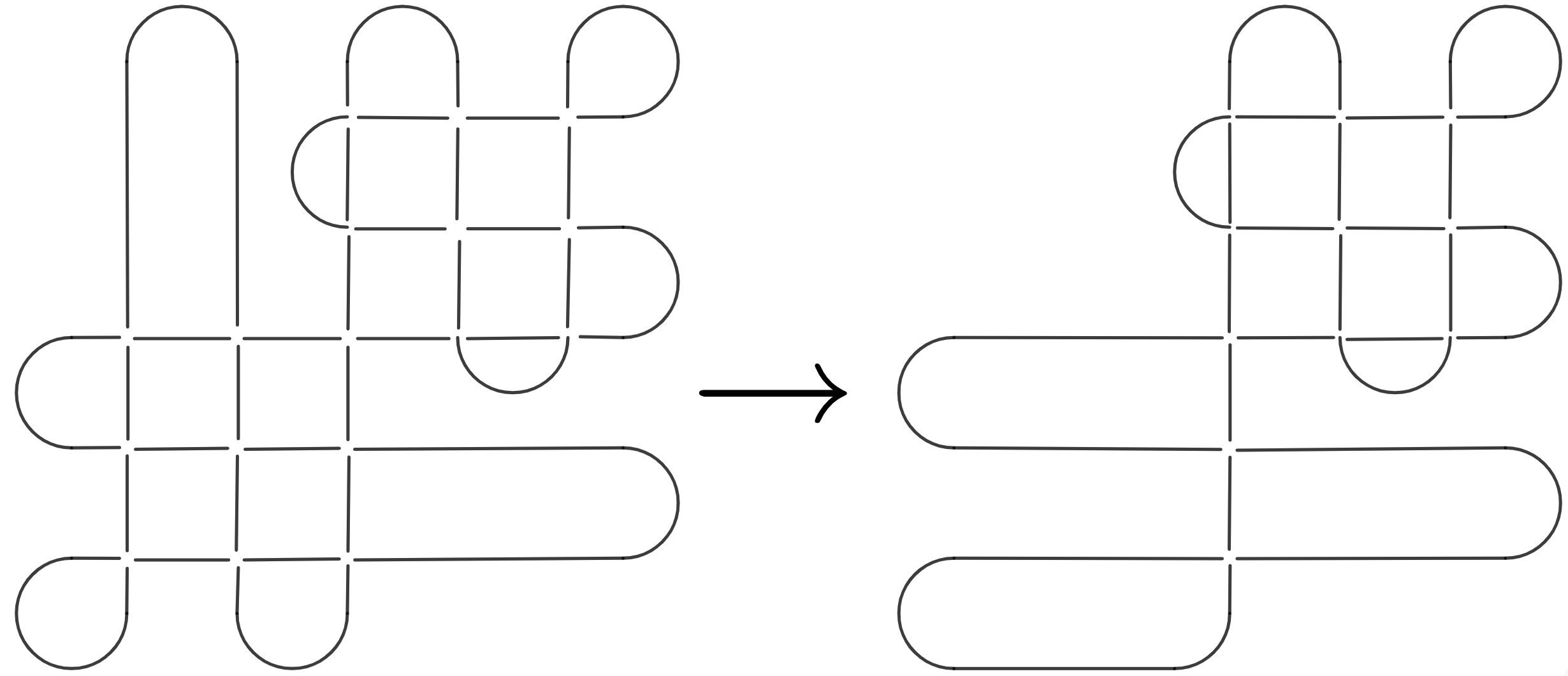}}
\end{minipage}
\hfill
\begin{minipage}[h]{0.5\linewidth}
\center{\includegraphics[width=0.9\linewidth]{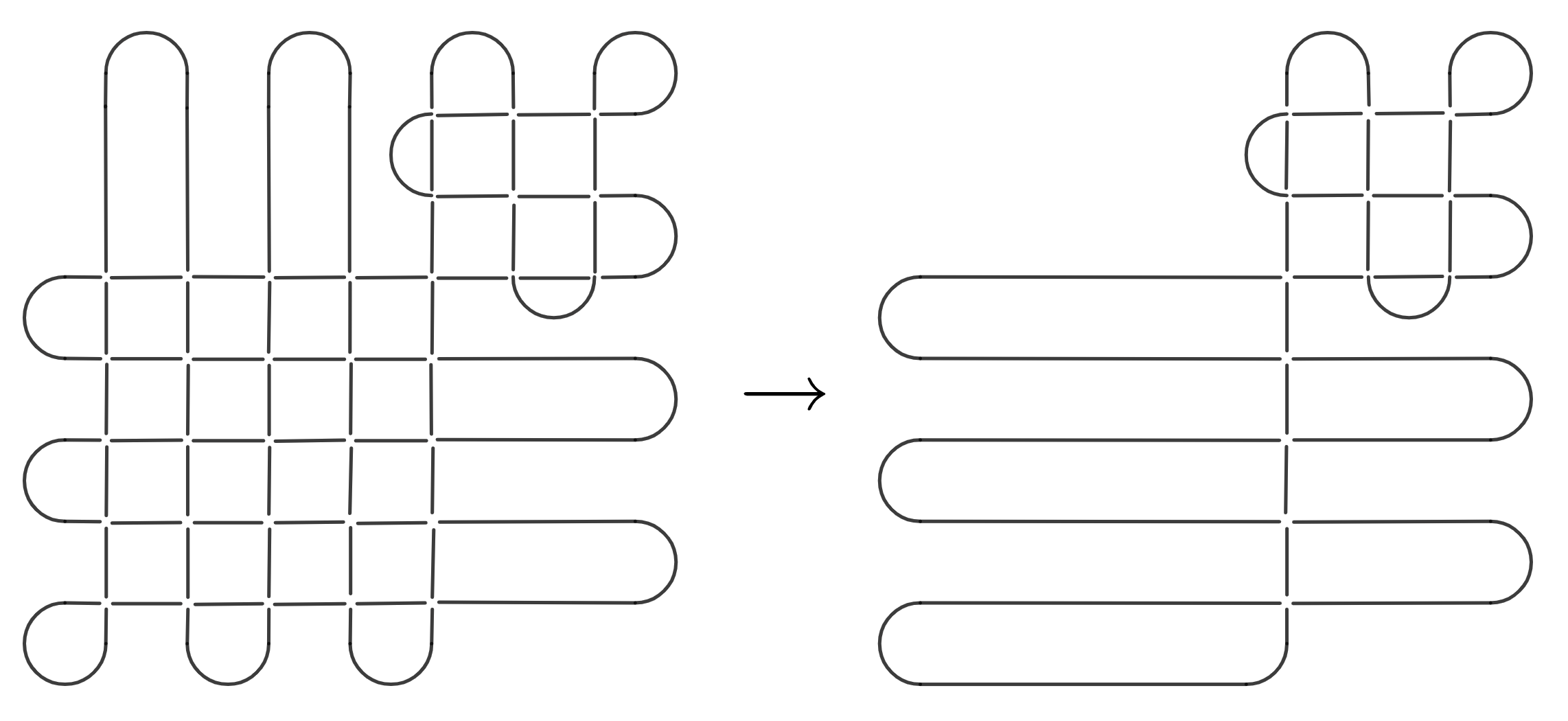}}
\end{minipage}
\caption{Vertical loops pull off. $3\times3$ smaller lattice is placed in the right upper corner; its intersections left empty as they can be chosen in order to form any lattice knot. Other intersections are to be determined in order to 'untie'\,.}
\end{figure}
    \item At the last step, we obtain a trivial loop, so the corresponding intersections can be arbitrary. Therefore, we get $2^{2(n-k)}$ abilities. 
\end{enumerate}
Actually, after point 1, one can untangle the remaining loop beginning with $n-k$ horizontal loops. This ability appears if $n-k>1$. We emphasise that one should exclude the variants obtained in the previous stages.
\begin{enumerate}
    \item[2'.] Pull off $n-k$ horizontal loops. This way, we get $2^{(n-k+2)(n-k)}$ other options. 
    \item[3'.] Finally, pull off $n-k$ vertical loops in $2^{2(n-k)}$ variants.
\end{enumerate}
Note that the position of the $(2k+1)\times(2k+1)$ lattice inside $(2n+1)\times(2n+1)$ lattice can be chosen $(n-k+1)\times(n-k+1)$ times, and the numbers of variants above stay the same even if the $(2k+1)\times(2k+1)$ lattice lies not on the boundary of $(2n+1)\times(2n+1)$ lattice (for example, see Fig.~\ref{NotOnBoundary}). Here it is crucial that $(2k+1)\times(2k+1)$ lattice is untiable because in the opposite case, we will overcount abilities to transform the bigger lattice knot into a smaller one.  
\\
In this analysis, we do not consider cases when the resulted smaller lattice knot is split into several parts separated in the bigger lattice. Therefore, we have the following\\
\textbf{Lemma 1.} There are totally more than $(n-k+1)^2[2^{2(n-k)(n+k+2)}+2^{(n-k)(n+3k+4)}(1-\delta_{n-k,1}-\delta_{n-k,0})]$ $(2k+1)\times(2k+1)$ lattices inside the $(2n+1)\times(2n+1)$ lattice.

\subsection{Knots inside lattices}\label{KnotsLat}
In order to find out how many knots with a fixed crossing number can be embedded into a lattice of size $(2n+1)\times(2n+1)$, one needs to find out which crossing numbers are allowed in a lattice with a fixed size.    

\begin{figure}[!ht]
    \centering
    \includegraphics[scale=0.43]{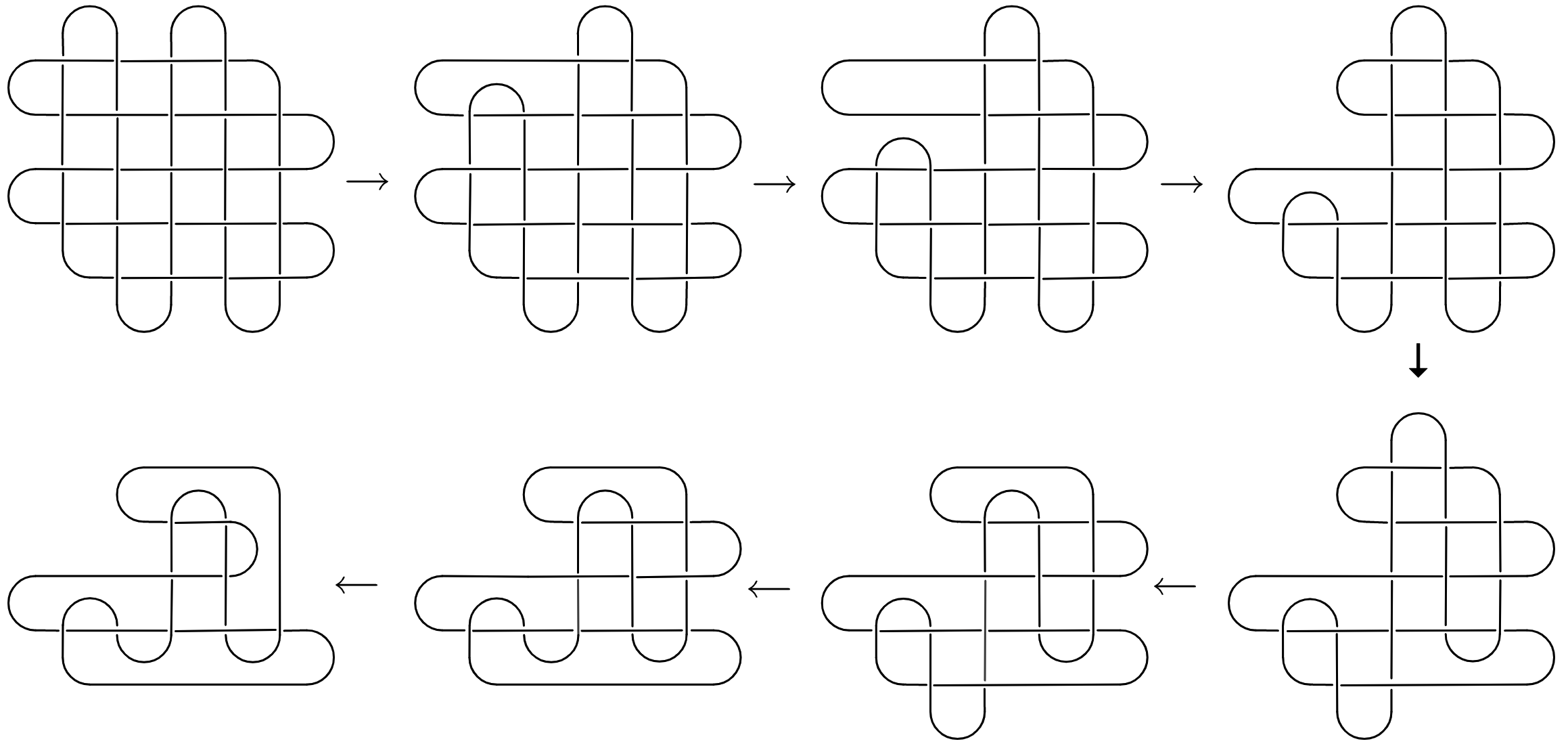}
    \caption{Knots with crossing numbers 23, 21, 19, 17, 15, 13, 11, and 9 inside the $5\times 5$ lattice.}\label{Ex}
\end{figure}

One can prove that there are knots with odd crossing numbers $\leq (2k+1)^2-2$ in the $(2k+1)\times(2k+1)$ lattice. Namely, one can easily find the biggest knot as the endless knot with a crossing number $(2k+1)^2-2$. And in order to obtain other knots, one pulls off loops of the endless knot, leaving the remaining crossings alternating. In Fig.~\ref{Ex}, there is an example (each step deletes two crossings, leaving the remaining alternating knot).
Note that there is no need to obtain other knots with smaller crossing numbers due to the theorem that all knots realized in a smaller lattice can be also realized in a bigger one. 
 \\ \\
Thus, due to the fact that one just needs to embed an untiable lattice of size $(2k+1)\times(2k+1)$ into the $(2n+1)\times(2n+1)$ lattice, \\
\textbf{Lemma 2.} the number of knots with $(2k-1)^2,(2k-1)^2+2,\dots,(2k+1)^2-4,(2k+1)^2-2$ crossing numbers (for $k=1$, the knot with crossing number 1 is the unknot and it must be discarded) in the $(2n+1)\times(2n+1)$ lattice (denote each of them as $\sigma_k^n$) can be estimated from below by the found number $(n-k+1)^2[2^{2(n-k)(n+k+2)}+2^{(n-k)(n+3k+4)}(1-\delta_{n-k,1}-\delta_{n-k,0})]$.

Let us obtain now the upper bound. As we have established, the amount of each of the knots with odd crossing numbers from $(2k-1)^2$ to $(2k+1)^2-2$ on a $(2n+1)\times(2n+1)$ lattice
\begin{equation}\label{smin}
    \sigma_k^n\geq (\sigma_k^n)_{\min} = (n-k+1)^2[2^{2(n-k)(n+k+2)}+2^{(n-k)(n+3k+4)}(1-\delta_{n-k,1}-\delta_{n-k,0})]\,.
\end{equation}
Number of different types of knots inside untiable $(2k+1)\times(2k+1)$ lattice is at least $8k$ for $k>1$, $6$ for $k=1$ and $1$ for $k=0$, so their total amount is
\begin{equation}\label{Smin}
\begin{cases}
    \Sigma_k^n\geq (\Sigma_k^n)_{\min} = 8k(n-k+1)^2[2^{2(n-k)(n+k+2)}+2^{(n-k)(n+3k+4)}(1-\delta_{n-k,1}-\delta_{n-k,0})]\,,\quad k>1\,, \\
    \Sigma_1^n\geq (\Sigma_1^n)_{\min} = 6n^2[2^{2(n-1)(n+3)}+2^{(n-1)(n+7)}(1-\delta_{n,2}-\delta_{n,1})]\,, \\
    \Sigma_0^n\geq (\Sigma_0^n)_{\min} = (\sigma_0^n)_{\min} = (n+1)^2[2^{2n(n+2)}+2^{n(n+4)}(1-\delta_{n,1}-\delta_{n,0})]\,.
\end{cases}
\end{equation}
Note that in the cases $k>0$, we take into account mirror knots, so that the amount of knots with fixed crossing numbers doubles. 

In all, there are $2^{(2n+1)^2}$ knots in the $(2n+1)\times(2n+1)$ lattice, as we can resolve each intersection in any of two possible variants. So, in order to obtain an upper bound on the amount of knots with a fixed crossing number, we need from $2^{(2n+1)^2}$ all knots subtract minimal number of knots with the rest crossing numbers, so that   
\begin{equation}
    \sigma_k^n\leq 2^{(2n+1)^2}- \sum_{i=0}^n (\Sigma_i^n)_{\min}\,+(\sigma_k^n)_{\min}.
\end{equation}
To sum up, we obtain the following \\
\textbf{Theorem 4.} The number of knots with $(2k-1)^2,(2k-1)^2+2,\dots,(2k+1)^2-4,(2k+1)^2-2$ crossing numbers (for $k=1$ the knot with crossing number 1 is the unknot and it must be discarded) in the $(2n+1)\times(2n+1)$ lattice satisfy the following constraints:
\begin{equation}\label{eval}
    (n-k+1)^2[2^{2(n-k)(n+k+2)}+2^{(n-k)(n+3k+4)}(1-\delta_{n-k,1}-\delta_{n-k,0})]\leq \sigma_k^n\leq 2^{(2n+1)^2}- \sum_{i=0}^n (\Sigma_i^n)_{\min}\,+(\sigma_k^n)_{\min}\,,
\end{equation}
where the minimal values are defined in~\eqref{smin} and~\eqref{Smin}. We emphasise that in this analysis, we differentiate between knots projections and their mirror ones.

\textbf{Remark 2.} According to our analysis, if there are knots with even crossing numbers $(2k-1)^2+1,(2k-1)^2+3,\dots,(2k+1)^2-3$ (for $k=1$ the knot with crossing number 2 is the unknot and it must be discarded), their numbers also admit the above bounds~\eqref{eval}.

Provide some examples:
\begin{equation}\label{ExEst}
\begin{aligned}
1 \leq\; &\sigma_1^1\leq 252\,, \\
256 \leq\; &\sigma_0^1\leq 506\,; \\
1 \leq\; &\sigma_2^2\leq 2^{25}-632\,847\,, \\
1024 \leq\; &\sigma_1^2\leq 2^{25}-631\,824\,, \\
626\,688 \leq\; &\sigma_0^2\leq 2^{25}-6144\,.
\end{aligned}
\end{equation}
Here $\sigma_0^1$ is the number of unknots inside $3\times 3$ lattice, and $\sigma_1^1$ is the amount of other types of knots inside $3\times 3$ lattice. Inside $5\times 5$ lattice, $\sigma_0^1$ is the number of unknots, $\sigma_1^2$ is the number of knots with crossing numbers $3,5,7$ and $4,6$ (if they appear), $\sigma_1^2$ is the number of knots with crossing numbers $9,11,13,15,17,19,21,23$ and $10,12,14,16,18,20,22$ (if they appear). From evaluation~\eqref{eval}, one can also see that unknots dominate for any $n$.

\setcounter{equation}{0}
\section{Numerical computations}\label{NumComp}
In order to check our estimates \eqref{eval}, we compute the Jones polynomials for small lattice knots. The Jones polynomials are good tools to distinguish small knots. The first example of knots with the identical Jones polynomials is $J^{5_1}(q)=J^{10_{132}}(q)$, and this becomes essential for knots with 10 and more crossings. 


\subsection{Knots inside 3x3 and 5x5 lattices}\label{3x35x5}

We have written a computer program which calculates the Jones polynomial of lattice knots using the state sum formula (see Subsection~\ref{Jones}). First, calculate which types of knots $3\times3$ lattice contains and their amounts. The computer program give us 10 different Jones polynomials. Comparing them with the Jones polynomials for knots with crossing numbers up to 7 (with the use of~\cite{knotinfo}) and using property~\eqref{JonesMir}, we identify:
\begin{equation}
\begin{aligned}
J^{0_1}(q)&=1\,, \\
J^{3_1}(q)&=J^{3_1^{mir}}(q^{-1})=q^2+q^6-q^8\,, \\
J^{4_1}(q)&=1 -q^2 - q^{-2} + q^4 + q^{-4}\,, \\
J^{5_2}(q)&=J^{5_2^{mir}}(q^{-1})=q^2-q^4+2 q^6-q^8+q^{10}-q^{12}\,, \\
J^{6_1}(q)&=J^{6_1^{mir}}(q^{-1})=q^{-4}-q^{-2}+2-2q^2+q^4-q^6+q^8\,, \\
J^{7_4}(q)&=J^{7_4^{mir}}(q^{-1})=q^2-2q^4+3 q^6-2q^8+3q^{10}-2q^{12}+q^{14}-q^{16}\,.
\end{aligned}
\end{equation}
So that we get the resulting Table~\ref{tab:3x3}. Note that the amounts of knots are in accordance with our estimates~\eqref{ExEst}. There are $2^9=512$ knots in $3\times3$ lattice. The ratio of unknots is $0.609$. Torus knots dominate, their ratio is $0.766$.
\begin{align}
    \centering
    \begin{tabular}{|c||c|c|c|c|c|c|}
    \hline
         knot & $0_1$ & $3_1$ & $4_1$ & $5_2$ & $6_1$ & $7_4$ \\
    \hline
         type & T & T & H & H & H & H \\
    \hline
         amount & 312 & 40+40 & 48 & 24+24 & 8+8 & 4+4 \\
    \hline
    \end{tabular}
    \label{tab:3x3}
\end{align}
   \noindent Table~\ref{tab:3x3}. Number of knots in $3\times 3$ lattice. In the 'type' raw, 'H' means hyperbolic knot and 'T' means torus one. If an amount of knots is split into a sum of two equal numbers, the corresponding knot is not amphichiral, so that there are a knot and its mirror one.

\bigskip

There are much more different knots inside $5\times 5$ lattice. We have computed $13\,829$ different Jones polynomials (up to $q\rightarrow q^{-1}$ change), thus there are not less then $13\,829$ different knots inside $5\times 5$ lattice. In Table~\ref{tab:5x5} we write down types and amounts of knots with crossing number less or equal 9. Note that the amounts of knots are in accordance with our estimates~\eqref{ExEst}. 

We have managed to identify only knots with up to 12 crossings as the known databases~\cite{katlas,knotinfo} contain the Jones polynomials for knots with crossing number not higher than 12. Thus, we are not able to find out the ratio of hyperbolic and satellite knots. The ratio of unknots is $0.151$, and they dominate. \\ 
\begin{align}
    \footnotesize{\centering
    \begin{tabular}{|c||c|c|c|c|c|c|c|c|c|c|c|c|c|c|c|}
    \hline
         knot & $0_1$ & $3_1$ & $4_1$ & $5_1$ & $5_2$ & $6_1$ & $6_2$ & $6_3$ &
         $7_1$ & $7_2$ & $7_3$ \\
    \hline
         type & T & T & H & T & H & H & H & H & T & H & H \\
    \hline
         amount & 5\,063\,616 & 2\,785\,728 & 1\,896\,896 & 525\,360 & 2\,327\,776 & 1\,253\,216 & 719\,496 & 381\,312 & 26\,560 & 508\,000 & 350\,784 \\
    \hline
    \hline
    knot & $7_4$ & $7_5$ & $7_6$ & $7_7$ & $8_1$ & $8_2$ & $8_3$ & $8_4$ & $8_5$ & $8_6$ & $8_7$ \\
    \hline
    type & H & H & H & H & H & H & H & H & H & H & H \\
    \hline
    amount & 519\,312 & 426\,256 & 475\,232 & 281\,280 & 250\,656 & 29\,824 & 130\,208 & 171\,648 & 26\,880 & 194\,816 & 42\,944 \\
    \hline
    \hline
    knot & $8_8$ & $8_9$ & $8_{10}$ & $8_{11}$ & $8_{12}$ & $8_{13}$ & $8_{14}$ & $8_{15}$
    & $8_{16}$ & $8_{17}$ & $8_{18}$ \\
    \hline
    type & H & H & H & H & H & H & H & H & H & H & H \\
    \hline
    amount & 249\,840 & 81\,536 & 49\,536 & 261\,184 & 130\,144 & 288\,704 & 332\,752 & 27\,456 & 37\,312 & 29\,440 & 13\,392 \\
    \hline
    \hline
    knot & $8_{19}$ & $8_{20}$ & $8_{21}$ & $9_1$ & $9_2$ & $9_3$ & $9_4$ & $9_5$ &
    $9_6$ & $9_7$ & $9_8$ \\
    \hline
    type & T & H & H & T & H & H & H & H & H & H & H \\
    \hline
    amount & 31\,680 & 69\,712 & 48\,064 & 0 & 89\,568 & 16\,192 & 21\,888 & 179\,072 & 13\,632 & 57\,088 & 101\,104 \\
    \hline
    \hline
    knot & $9_9$ & $9_{10}$ & $9_{11}$ & $9_{12}$ & $9_{13}$ & $9_{14}$ & $9_{15}$ & $9_{16}$ & $9_{17}$ & $9_{18}$ & $9_{19}$ \\
    \hline
    type & H & H & H & H & H & H & H & H & H & H & H \\
    \hline
    amount & 27\,968 & 65\,392 & 16\,000 & 169\,152 & 157\,568 & 117\,584 & 111\,936 & 14\,112 & 8\,128 & 153\,072 & 142\,864 \\
    \hline
    \hline
    knot & $9_{20}$ & $9_{21}$ & $9_{22}$ & $9_{23}$ & $9_{24}$ & $9_{25}$ & $9_{26}$ & $9_{27}$ & $9_{28}$ & $9_{29}$ & $9_{30}$ \\
    \hline
    type & H & H & H & H & H & H & H & H & H & H & H \\
    \hline
    amount & 35\,008 & 165\,568 & 44\,352 & 98\,208 & 28\,800 & 25\,344 & 32\,896 & 40\,128 & 15\,552 & 3\,680 & 42\,160 \\
    \hline
    \hline
    knot & $9_{31}$ & $9_{32}$ & $9_{33}$ & $9_{34}$ & $9_{35}$ & $9_{36}$ & $9_{37}$ & $9_{38}$ & $9_{39}$ & $9_{40}$ & $9_{41}$ \\
    \hline
    type & H & H & H & H & H & H & H & H & H & H & H \\
    \hline
    amount & 17\,312 & 32\,256 & 16\,320 & 14\,432 & 9\,504 & 39\,104 & 19\,200 & 3\,792 & 16\,304 & 336 & 5\,728 \\
    \hline
    \hline
    knot & $9_{42}$ & $9_{43}$ & $9_{44}$ & $9_{45}$ & $9_{46}$ & $9_{47}$ & $9_{48}$ & $9_{49}$ & & & \\
    \hline
    type & H & H & H & H & H & H & H & H & & & \\
    \hline
    amount & 44\,928 & 57\,312 & 64\,992 & 50\,464 & 35\,712 & 12\,608 & 15\,592 & 7\,264 & & & \\
    \hline
    \end{tabular}}
    \label{tab:5x5}
\end{align}
\noindent Table~\ref{tab:5x5}. Number of knots in $5\times 5$ lattice with up to 9 crossings. To shorten the table, we do not differ between knots and their mirrors.

\bigskip

\begin{figure}[h]
\begin{minipage}[h]{0.5\linewidth}
\center{\includegraphics[width=0.95\linewidth]{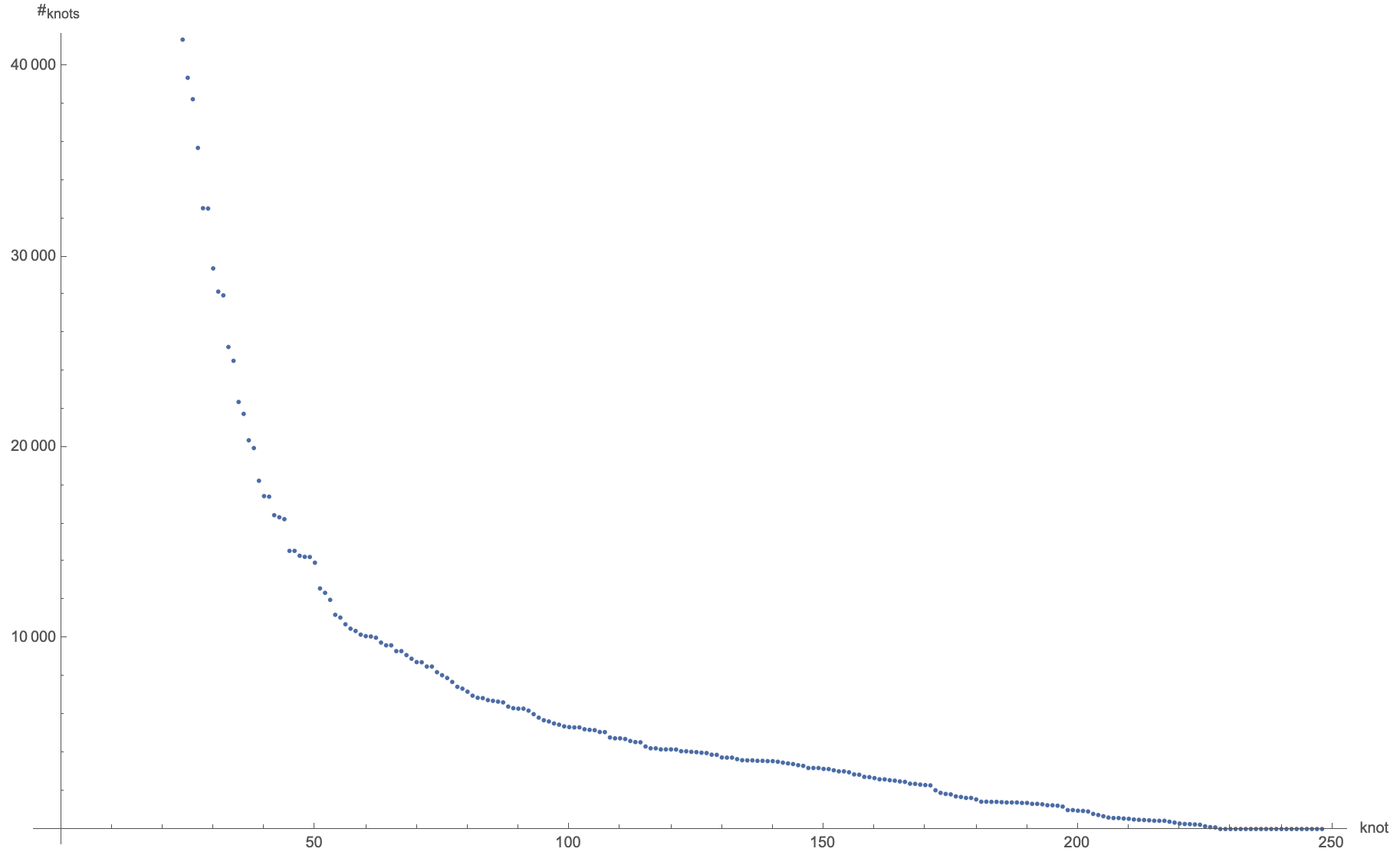}}
\end{minipage}
\hfill
\begin{minipage}[h]{0.5\linewidth}
\center{\includegraphics[width=0.95\linewidth]{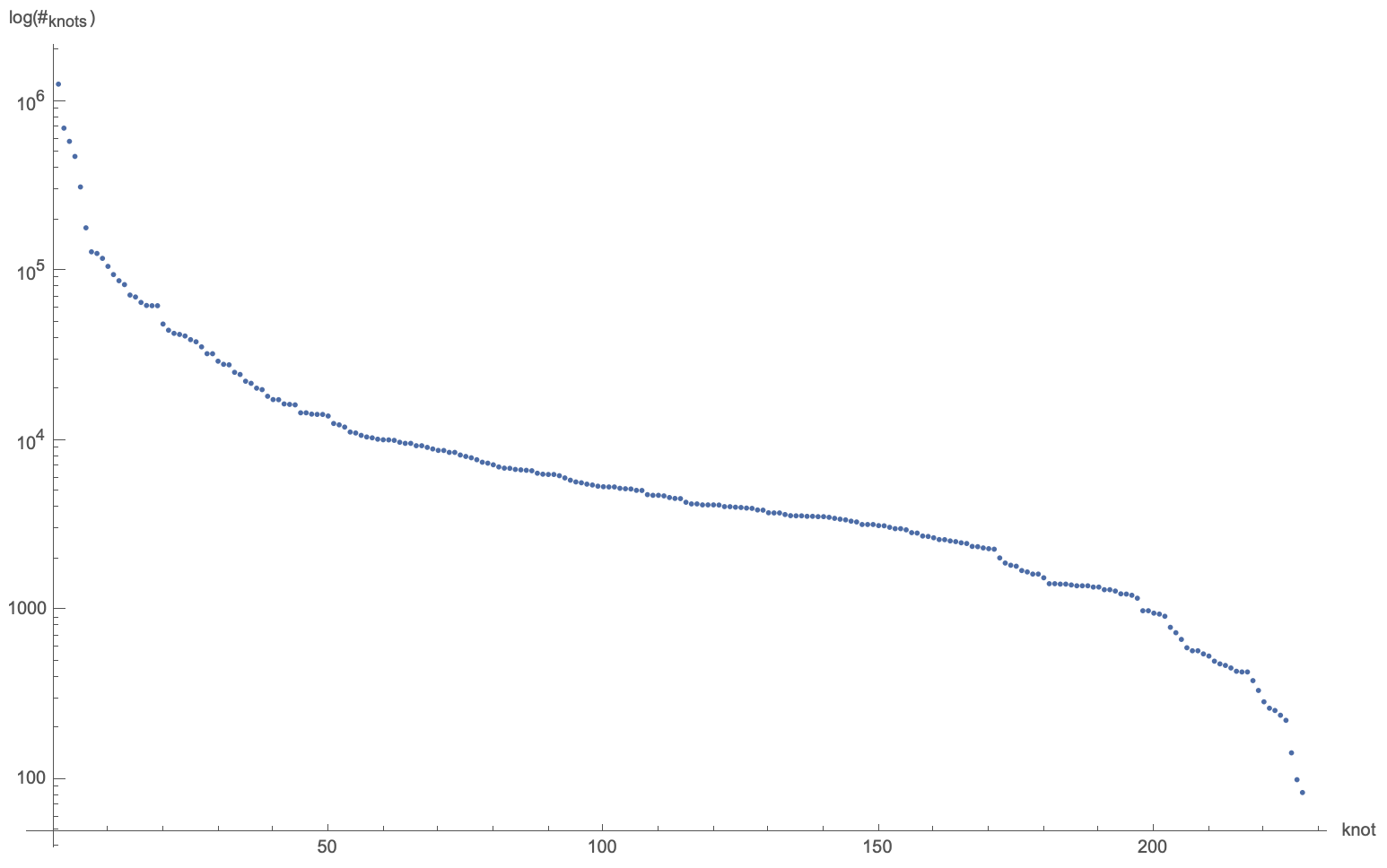}}
\end{minipage}
\caption{The amounts of knots with up to 10 crossings in $5\times 5$ lattice, ordinary and logarithmic scales. The horizontal axis just enumerates knots.}\label{plots}
\end{figure}

Note that the amounts of knots in $5\times 5$ lattice decrease almost exponentially, Fig.~\ref{plots}.

{\subsection{Ratio of unknots} \label{sec:ratio-of-unknots}
    Above we have provided results of extensive numeric experiments on $3 \times 3$
    and $5 \times 5$ lattices. To go to lattices of bigger size one needs to do two adjustments
(otherwise the numeric computations required become overwhelming even for modern computers).
\begin{itemize}
\item Limit the scope: switch from detecting all knots to just detecting the unknot, and
whether the knot is outside the Rolfsen table (i.e. has more than 10 minimal crossings). As mentioned above,
this is done with help of evaluating the Jones polynomial so, due to unavoidable collisions the result
will be estimate from above. This optimization drastically reduces memory requirements as we
don't keep track of all different knots.
\item Reduce precision: instead of honestly iterating over all knots on a given lattice,
generate a (reasonably) large number of random knots on this lattice (\textit{a la} Monte-Carlo).
From this sampling one gets an estimate of the ratio of unknots, rather than accurate total number.
With this optimization we can set computation time to any (reasonable) desired value.
\end{itemize}

With these optimizations we get the following estimates for the number of unknots on a given lattice: 
  \begin{align}
  \begin{array}{|c|c|c|}
  \hline \text{Lattice size} & \text{Ratio of unknots} & \text{Ratio of $> 10$ crossings} \\
  \hline 1 \times 1 & 1.0 & 0.0 \\
  \hline 3 \times 3 & 0.614 & 0.0 \\
  \hline 5 \times 5 & 0.134 & 0.272 \\
  \hline 7 \times 7 & 0.015 & 0.872 \\
  \hline 9 \times 9 & 0.0001 & 0.9972 \\
  \hline
  \end{array}
  \end{align}
    We see that: 
    \begin{itemize}
    \item    Firstly, these numbers are in accordance (modulo imprecisions introduced
    by random sampling and Jones collisions) with theoretical estimates
    \eqref{eval} and full-iteration numerics (see Tables~\ref{tab:3x3} and \ref{tab:5x5}),
    \item Secondly, the ratio of unknots
    seems to be exponentially falling, confirming to (asymptotic mean-field) estimates
    of Grosberg and Nechaev~\cite{Grosberg_1992,nechaev1999statistics}.
    \item Thirdly, the number of non-Rolfsen (i.e. with $> 10$ minimal crossings) knots
    in a rectangular grid is rapidly (exponentially) approaching 1, that is, the task of finding
    a knot with small number of crossings on a large grid is exponentially complex
    computational task. This is also in accord with our theoretical estimates \eqref{eval}.
    The implications of this will be explored in the future research.
    \end{itemize}

}

\setcounter{equation}{0}
\section{Conclusion/Discussion}\label{conclusion}
In this paper, we have answered classical questions of knot theory regarding the introduced stratification by lattice knots. Namely, we have obtained the following results.
\begin{enumerate}
    \item We have got evaluations on the amounts of knots with fixed crossing number inside $(2n+1)\times(2n+1)$ lattice $\forall\, n$ (Section~\ref{Bounds}).
    \item We have classified (according to the Rolfsen table) all knots in $3\times 3$ lattice (Table~\ref{tab:3x3}) and knots up to 12 crossings which can be distinguished by the Jones polynomial in $5\times 5$ lattice (there are only knots with up to 9 crossings in Table~\ref{tab:5x5}), Section~\ref{3x35x5}.
    \item We have found amounts of unknots for five lattices (Section~\ref{sec:ratio-of-unknots}). Their ratio decreases exponentially with the growth of a lattice size.
    \item We have found amounts of knots with $>10$ minimal crossings for five lattices (Section~\ref{sec:ratio-of-unknots}). Their ratio increases exponentially with the growth of a lattice size.
\end{enumerate}
In order to find out which type of knots dominate for big lattice size, we need to study how one can effectively calculate knot invariants which distinguish between hyperbolic, torus and satellite knots.

With the growth of lattice size there appear more knots which cannot be distinguished by the Jones polynomial. Moreover, for the Jones polynomial, it is even not proved whether it distinguish the unknot. It is more effective to use the Khovanov polynomial~\cite{Bar_Natan_2002,Kh1}, that is a categorification of the Jones polynomial, as there is a theorem that the Khovanov polynomial detects the unknot~\cite{kronheimer2011khovanov}. Moreover, it distinguishes more knots. However, it is much more complicated to calculate the Khovanov polynomial, and the development of simpler method for its calculation for lattice knots is a separate creative problem.

\setcounter{equation}{0}
\section*{Acknowledgements}\label{Acknowledgements}
We would like to thank A. Malyutin and Yu. Belousov for useful discussions. This work was funded by the grant of Leonhard Euler International Mathematical Institute in Saint Petersburg № 075–15–2019–1619 (E.L., N.T.), by the grants of the Foundation for the Advancement of Theoretical Physics and Mathematics “BASIS” (E.L., N.T.), by the RFBR grant 20-01-00644 (N.T., A.P.), by the joint RFBR and TUBITAK grant 21-51-46010-CT\_a (N.T.) and by the joint RFBR and MOST grant 21-52-52004\_MHT (A.P.).

\printbibliography

\end{document}